\documentclass[12pt,reqno]{amsart}
\usepackage{amsmath}
\usepackage{amssymb}
\usepackage{graphicx}
\usepackage{amsfonts}
\usepackage[margin=0.93in, letterpaper ]{geometry}
\usepackage{lineno}
\usepackage[onehalfspacing]{setspace}

\newtheorem{theorem}{Theorem}
\theoremstyle{plain}

\newtheorem{claim}[theorem]{Claim}

\newtheorem{corollary}[theorem]{Corollary}

\newtheorem{problem}[theorem]{Problem}
\newtheorem{proposition}[theorem]{Proposition}

\newtheorem{remark}[theorem]{Remark}

\newtheorem{scholium}[theorem]{Scholium}
\numberwithin{equation}{section}
\numberwithin{theorem}{subsection}
\input{tcilatex}

\begin{document}
\title{The Pieri Rule for $GL_{n}$ Over Finite Fields }
\author{\textsf{Shamgar Gurevich}}
\address{\textit{Department of Mathematics, University of Wisconsin,
Madison, WI 53706, USA.}}
\email{shamgar@math.wisc.edu}
\author{\textsf{Roger Howe}}
\address{\textit{Department of Mathematics, Yale University, New Haven, CT
06520, USA.}}
\email{roger.howe@yale.edu}

\begin{abstract}
The Pieri rule gives an explicit formula for the decomposition of the tensor
product of irreducible representation of the complex general linear group $%
GL_{n}(%
\mathbb{C}
)$ with a symmetric power of the standard representation on $%
\mathbb{C}
^{n}$. It is an important and long understood special case of the
Littlewood-Richardson rule for decomposing general tensor products of
representations of $GL_{n}(%
\mathbb{C}
)$.

In our recent work \cite{Gurevich-Howe17, Gurevich-Howe19} on the
organization of representations of the general linear group over a finite
field $\mathbb{F}_{q}$ using small representations, we used a generalization
of the Pieri rule to the context of this latter group.

In this note, we demonstrate how to derive the Pieri rule for $GL_{n}(%
\mathbb{F}_{q})$. This is done in two steps; the first, reduces the task to
the case of the symmetric group $S_{n}$, using the natural relation between
the representations of $S_{n}$ and the \textit{spherical principal series }%
representations of $GL_{n}(\mathbb{F}_{q})$; while in the second step,
inspired by a remark of Nolan Wallach, the rule is obtained for $S_{n}$
invoking the $S_{l}$-$GL_{n}(%
\mathbb{C}
)$ Schur duality.

Along the way, we advertise an approach to the representation theory of the
symmetric group which emphasizes the central role played by the \textit{%
dominance order} on Young diagrams. The ideas leading to this approach seem
to appear first, without proofs, in \cite{Howe-Moy86}.
\end{abstract}

\maketitle
\dedicatory{\smallskip\ \ \ \ \ \ \ \ \ \ \ \ \ \ \ \ \ \ \ \ \ \ \ \ \ \ \
\ \ \ \textrm{Dedicated to the memory of Ronald Douglas}}

\section{\textbf{Introduction}}

Two basic tasks in the representation theory of a finite group $G$ are: the
parameterization of its set $\widehat{G}$ (of isomorphism classes) of
irreducible representations (irreps); and the decomposition into direct sum
of irreps of certain of its naturally arising representations.

The Pieri rule that we formulate and prove in this note addresses a
particular instance of the second task mentioned above, for the case of the
general linear group $GL_{n}=GL_{n}(\mathbb{F}_{q})$ over a finite field $%
\mathbb{F}_{q}$. It can be used to give a recursive solution to the general
problem of decomposing the permutation actions of $GL_{n}$ on functions on
flag manifolds.

The Pieri rule can be useful in other ways. Indeed, in \cite%
{Gurevich-Howe17, Gurevich-Howe19} we developed a precise notion of "size"
for irreps of $GL_{n},$ called "tensor rank". This is an integer $0\leq
k\leq n$ that is naturally attached to an irreducible representation (irrep)
and helps to compute important analytic properties such as its dimension and
character values on certain elements of interest. In particular, in \textit{%
loc. cit. }the Pieri rule for $GL_{n}$ enabled us to give an effective
formula for the irreps of $GL_{n}$ of a given tensor rank $k$.

We proceed to consider the subgroups involved in the construction of
representations involved in the formulation of the Pieri rule.

\subsection{\textbf{Young Diagrams and Parabolic Subgroups}}

The representations we are interested in are naturally realized on spaces
constructed using standard parabolic subgroups \cite{Borel69} of general
linear groups, as we will now describe.

Fix an integer $0\leq k\leq n,$ and denote by $\mathcal{Y}_{k}$ the
collection of Young diagrams of size $k$ \cite{Fulton97}. In more detail, by
a \textit{Young diagram }(or\textit{\ partition})\textit{\ }$D\in \mathcal{Y}%
_{k}$, we mean an ordered list of non-negative integers 
\begin{equation}
D=(d_{1}\geq ...\geq d_{r})\text{, with }d_{1}+...+d_{r}=k.  \label{D}
\end{equation}

\begin{figure}[h]\centering
\includegraphics
{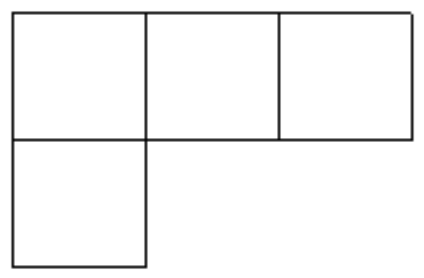}%
\caption{The Young diagram $D=(3,1)\in \mathcal{Y}_{4}$. }\label{young-d-3-1}%
\end{figure}

It is common to visualize---see Figure \ref{young-d-3-1} for
illustration---the diagram $D$ with the help of a drawing of $r$ rows of
square boxes, each row one on top of the other, starting at the left upper
corner, in such a way that the $i$-th row contains $d_{i}$ boxes.%

\begin{figure}[h]\centering
\includegraphics
{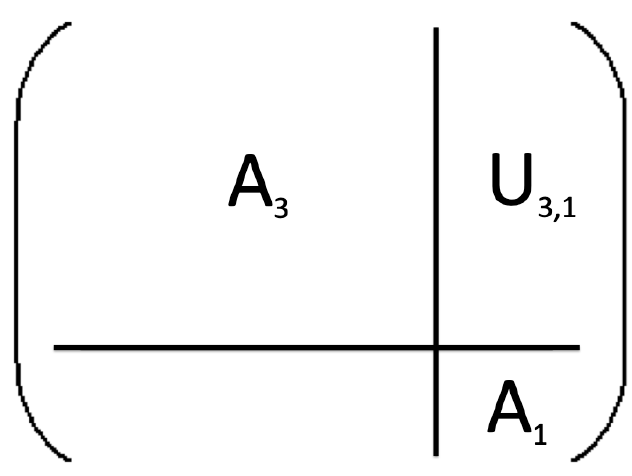}%
\caption{The parabolic $P_{D}\subset GL_{4}$, $D=(3,1)$, has $A_{3}\in
GL_{3},$ $A_{1}\in GL_{1}$, $U_{3,1}\in M_{3,1}$. }\label{pd31}%
\end{figure}

To the diagram $D$ (\ref{D}) we can attach the following increasing sequence 
$F_{D}$ of subspaces of the $k$-dimensional vector space $\mathbb{F}_{q}^{k}$%
:%
\begin{equation}
F_{D}:0\subset \mathbb{F}_{q}^{d_{1}}\subset \mathbb{F}_{q}^{d_{1}+d_{2}}%
\subset \ldots \subset \mathbb{F}_{q}^{k},  \label{FD}
\end{equation}%
and call it the \textit{standard flag }attached to $D$. In particular,
having $D$ we can form---see Figure \ref{pd31} for illustration\footnote{%
We denote $M_{k,n}=M_{k,n}(\mathbb{F)}$ the space of $k\times n$ matrices
over a field $\mathbb{F}$.}---the stabilizer subgroup%
\begin{equation}
P_{D}=Stab_{GL_{k}}(F_{D})\subset GL_{k},  \label{PD}
\end{equation}%
that we will call the \textit{standard parabolic subgroup attached to} $D$.

Probably the most important example from this class of parabolic subgroups
is the Borel subgroup $B$ of upper triangular matrices in $GL_{k}$, which is
just $P_{D}$ with%
\begin{equation*}
D=\left. 
\begin{tabular}{c}
\hline
\multicolumn{1}{|c|}{} \\ \hline
\multicolumn{1}{|c|}{} \\ \hline
$\vdots $ \\ 
$\vdots $ \\ \hline
\multicolumn{1}{|c|}{} \\ \hline
\end{tabular}%
\right\} \text{ }k\text{ times.}
\end{equation*}

Next, we describe the specific type of representations that the Pieri rule
attempts to decompose.

\subsection{\textbf{The Pieri Problem\label{S-TPP}}}

Take $D\in \mathcal{Y}_{k},$ denote by $\mathbf{1}$ the trivial
representation of $P_{D}$, and consider the induced representation 
\begin{equation*}
I_{D}=Ind_{P_{D}}^{GL_{k}}(\mathbf{1),}
\end{equation*}%
which is given by the space of complex valued functions on $GL_{k}/P_{D},$
equipped with the standard left action of $GL_{k}$ on it.

In the case when $P_{D}=B$ is the Borel subgroup, the set $GL_{k}/B$ is the
flag variety, and we will call the collection of irreps that appear inside $%
Ind_{B}^{GL_{k}}(\mathbf{1})$ the \textit{spherical principal series (SPS)}.

There is a natural recipe (that we will recall in detail below) to
parametrize the \textit{SPS} by Young diagrams. Note that for each Young
diagram $D\in \mathcal{Y}_{k},$ we have $I_{D}<Ind_{B}^{GL_{k}}(\mathbf{1})$%
, where $<$ denotes subrepresentation. Interestingly, each $I_{D}$ contains
(with multiplicity one) a well defined "largest" irreducible
subrepresentation $\rho _{D}$. We will leave the details of that story for
the body of the note, but the collection $\{\rho _{D};$ $D\in \mathcal{Y}%
_{k}\}$ realizes the totality of SPS representations of $GL_{k}$.

We proceed to formulate the Pieri problem.

Fix $0\leq k\leq n$, and denote by $P_{k,n-k}\subset GL_{n}$ the standard
parabolic fixing the first $k$ coordinate subspace of $\mathbb{F}_{q}^{n}$.
There is a natural surjective map $P_{k,n-k}\twoheadrightarrow GL_{k}\times
GL_{n-k}$. Take an SPS\ representation $\rho _{D}$ of $GL_{k}$, and denote
by $\mathbf{1}_{n-k}$ the trivial representation of $GL_{n-k}$. Pull back
the representation $\rho _{D}\otimes \mathbf{1}_{n-k}$ from $GL_{k}\times
GL_{n-k}$ to $P_{k,n-k}$ and form the induced representation 
\begin{equation}
I_{\rho _{D}}=Ind_{P_{k,n-k}}^{GL_{n}}(\rho _{D}\otimes \mathbf{1}_{n-k}).
\label{IrhoD}
\end{equation}

Now we can write down the natural,

\begin{problem}[Pieri problem]
\label{P-PP}Decompose the representation $I_{\rho _{D}}$ into irreducibles.
\end{problem}

It is easy to see that the components of $I_{\rho _{D}}$ (\ref{IrhoD}) are
SPS representations of $GL_{n}$, so we are looking for a solution to Problem %
\ref{P-PP} in terms of Young diagrams, i.e., members of $\mathcal{Y}_{n}$.

In this note we present a solution to the Pieri problem for $GL_{n}$ in two
steps. First we explain why it is enough to solve the analogous problem for
the representations of the symmetric group $S_{n}$. Then, in the second
step, we demonstrate that the Pieri rule holds for $S_{n}$, invoking the
Schur (a.k.a. Schur-Weyl) duality for $S_{l}$-$GL_{n}(%
\mathbb{C}
)$, and a use of the classical Pieri rule for $GL_{n}(%
\mathbb{C}
)$ \cite{Howe92, Pieri1893, Weyman89}. We note that in \cite%
{Ceccherini-Silberstein-Scarabotti-Tolli10}, Section 3.5, there is a proof
of the Pieri rule for $S_{n}$ based on a quite different approach.\smallskip

We proceed to a short acknowledgements paragraph, after that give the table
of contents, and start the body of the note.\medskip

\textbf{Acknowledgements. }The material presented in this note is based upon
work supported in part by the National Science Foundation under Grants No.
DMS-1804992 (S.G.) and DMS-1805004 (R.H.). \ 

This note was written during 2019 while S.G. was visiting the Math
Department and the College of Education at Texas A\&M University, and he
would like to thank these institutions, and to thank personally R. Howe.

\tableofcontents

\section{\textbf{Representations of }$S_{n}$\label{S-RepSn}}

The standard parametrization of the irreps of $S_{n}$ is done using Young
diagrams \cite{Sagan91}. We will discuss various aspects of the construction
leading to this parametrization, emphasizing the role played by the
dominance relation on the set $\mathcal{Y}_{n}$ of Young diagrams. We will
follow closely ideas formulated (without proofs) in Appendix 2 of \cite%
{Howe-Moy86}.

\subsection{\textbf{The Young Modules\label{S-YM}}}

Recall that partitioning the set $\{1,..,n\}$ into $r$ disjoint subsets of
size $d_{i}$ each, and assigning these numbers, respectively, to the rows of
the Young diagram $D=(d_{1}\geq ...\geq d_{r})\in \mathcal{Y}_{n},$ gives
rise to a \textit{Young tabloid }\cite{Fulton97}. Let us denote by $\mathcal{%
T}_{D}$ the collection of all Young tabloids that one can make using $D$.
The natural action of the group $S_{n}$ on $\mathcal{T}_{D}$ is transitive.
Moreover, we can identify 
\begin{equation*}
\mathcal{T}_{D}=S_{n}/S_{D},
\end{equation*}%
where $S_{D}\subset S_{n}$ is the stabilizer subgroup 
\begin{equation}
S_{D}=Stab_{S_{n}}(T_{D}),  \label{SD}
\end{equation}%
of the tabloid $T_{D}$ that obtained by assigning to the first row of $D$
the numbers $1,..,d_{1},$ to the second $d_{1}+1,...,d_{1}+d_{2}$, etc. The
group $S_{D}$ is naturally isomorphic to the product $S_{d_{1}}\times
...\times S_{d_{r}}$ embedded in $S_{n}$ in the usual way.

Now, we consider the induced representation, called the \textit{Young module
associated to }$D,$ 
\begin{equation}
Y_{D}=Ind_{S_{D}}^{S_{n}}(\mathbf{1}),  \label{YD}
\end{equation}%
where $\mathbf{1}$ stands for the trivial representation of $S_{D}$. It is
naturally realized as the permutation representation of $S_{n}$ on the space
of functions on $\mathcal{T}_{D}$.

\subsection{\textbf{Properties of the Young Modules}}

We derive basic properties of the family of Young modules (\ref{YD}). They
give, in particular, as a corollary the standard classification of the
irreps of $S_{n}$, and, as we mentioned earlier, they can be effectively
understood using the important \underline{dominance} relation $\preceq $ on
the set $\mathcal{Y}_{n}$ of Young diagrams, which we recall now.

\begin{figure}[h]\centering
\includegraphics
{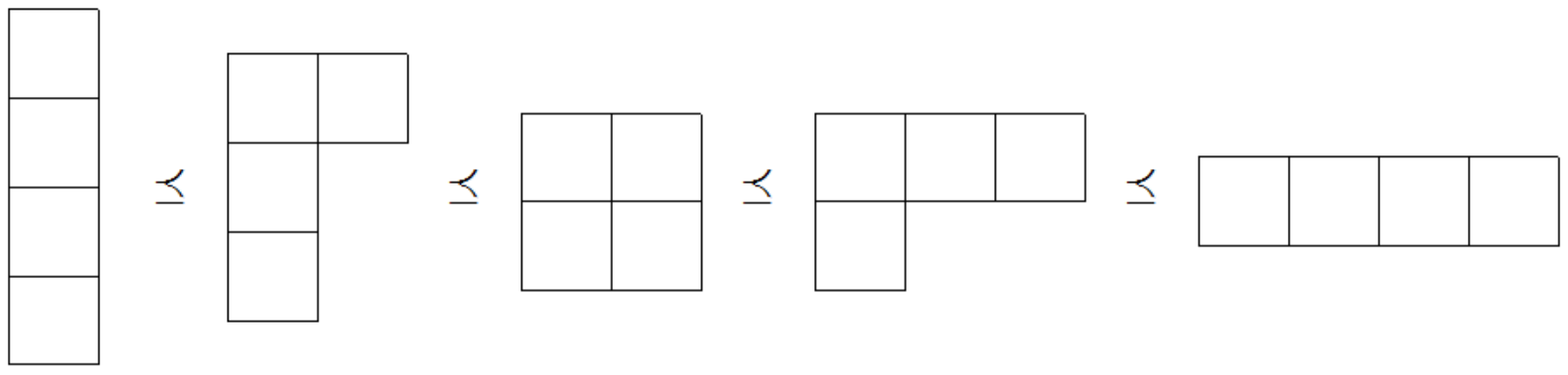}%
\caption{The set $\mathcal{Y}_{4}$ is totally ordered by $\preceq $ (this is
not true for $\mathcal{Y}_{n},$ $n\geq 6$).}\label{y4}%
\end{figure}

Suppose---see Figure \ref{y4} for illustration---we have a Young diagram $D$
which is obtained from another diagram by moving one of the boxes of $%
D^{\prime }$ to a (perhaps new) lower row, then we write 
\begin{equation}
D\preceq D^{\prime },  \label{DR}
\end{equation}%
and $\preceq $ on $\mathcal{Y}_{n}$ is the order generated from all the
inequalities of the form (\ref{DR}).%

Now, using the terminology afforded by the dominance relation, we can
formulate the main technical results concerning the Young modules.

For two representations $\pi $ and $\tau $ of a finite group $G$, let us
denote by $\left\langle \pi ,\tau \right\rangle $ their \textit{intertwining
number \cite{Serre77} }%
\begin{equation}
\left\langle \pi ,\tau \right\rangle =\dim Hom(\pi ,\tau ).  \label{IN}
\end{equation}
In addition, we denote the \textit{sign }representation of $S_{n}$ by $sgn$,
and introduce the \textit{twisted} Young module $%
Y_{E}(sgn)=Ind_{S_{E}}^{S_{n}}(sgn)$ attached to $E\in \mathcal{Y}_{n}.$
Finally, let us denote by $D^{t}$ the diagram in $\mathcal{Y}_{n}$ which is 
\textit{transpose} to $D$. That is, $D^{t}$ is gotten from $D$ by reflecting
across the downward diagonal from the top left box; in other words, the
columns of $D$ become the rows of $D^{t}$. Then,

\begin{theorem}
\label{T-Main}For any two Young diagrams $D,E\in \mathcal{Y}_{n}$, we have,

\begin{enumerate}
\item \textit{Intertwinity}: $\ \left\langle Y_{E}(sgn),Y_{D}\right\rangle
=\left\{ 
\begin{array}{c}
0,\text{ iff }E\npreceq D^{t}; \\ 
1,\text{ if \ }E=D^{t},%
\end{array}%
\right. $

and,

\item \textit{Monotonicity:}$\ \ D\precneqq E$ if and only if $Y_{E}\lneqq
Y_{D}$.
\end{enumerate}
\end{theorem}

For a proof of Theorem \ref{T-Main}, see Appendix \ref{P-T-Main}.\smallskip\ 

\subsection{\textbf{The Irreducible Representations of }$S_{n}$}

Part (1) of Theorem \ref{T-Main} produces the standard classification, by
Young diagrams, of the unitary dual (i.e., the set of irreps) $\widehat{S}%
_{n}$ of $S_{n},$ due to Frobenius and others \cite{Frobenius68}. Indeed,
for each $D\in \mathcal{Y}_{n},$ let us denote by $\sigma _{D}$ the unique
joint component of $Y_{D^{t}}(sgn)$ and $Y_{D}$. Then,

\begin{corollary}[Classification]
\label{C-Class}The irreps 
\begin{equation}
\sigma _{D},\text{ \ }D\in \mathcal{Y}_{n},  \label{sigmaD}
\end{equation}%
are pairwise non-isomorphic and exhaust $\widehat{S}_{n}$.
\end{corollary}

For a proof of Corollary \ref{C-Class} see \ref{P-C-Class}.

\subsection{\textbf{The Grothendieck Group of }$S_{n}$}

In Section \ref{S-SPS} we will draw certain conclusions for the
representation theory of the general linear group $GL_{n}=GL_{n}(\mathbb{F}%
_{q})$, using the properties obtained in this section for the
representations of $S_{n}$. An effective way to formulate this passage from $%
S_{n}$ to $GL_{n},$ is to use the formalism of the Grothendieck group of
representations, and in particular to describe consequences of Theorem \ref%
{T-Main} to the structure of this group in the case of $S_{n}$.

Given a finite group $G$, we can consider the Abelian group $K(G)$ generated
from the set $\widehat{G}$ of isomorphism classes of irreps of $G$ using the
direct sum operation $\oplus $. Note that $K(G)$ has a natural partial order 
$<$ given by the sub-representation relation, and it comes equipped with a
bilinear form $\left\langle \text{ },\right\rangle ,$ giving any two
representations $\pi ,\tau $, their intertwining number $\left\langle \pi
,\tau \right\rangle $ (\ref{IN}).

In particular, $K(S_{n})$ is a free $%
\mathbb{Z}
$-module with basis $\widehat{S_{n}}=\{\sigma _{D},$ \ $D\in \mathcal{Y}%
_{n}\}$, where $\sigma _{D}$ are the irreps (\ref{sigmaD}). \ However, $%
K(S_{n})$ has another natural $%
\mathbb{Z}
$-basis, i.e.,

\begin{proposition}
\label{P-Basis}The collection of Young modules $Y_{D},$ $D\in \mathcal{Y}%
_{n},$ forms a $%
\mathbb{Z}
$-basis for $K(S_{n})$.
\end{proposition}

Proposition \ref{P-Basis} follows from the following two consequences of
Theorem \ref{T-Main}:

\begin{scholium}
\label{Sc-SC}The following hold,

\begin{enumerate}
\item \textit{Spectrum: }The irrep $\sigma _{E}$ (\ref{sigmaD}), appears in
the Young module $Y_{D}$ if and only if $D\preceq E$.

\item \textit{Characterization:} The irrep $\sigma _{D}$ (\ref{sigmaD}) is
the only irrep that appears in $Y_{D}$ but not in $Y_{E}$ for every $%
D\precneqq E$.
\end{enumerate}
\end{scholium}

In particular, from Part (2) of Scholium \ref{Sc-SC} we deduce that the
collection of Young modules is a minimal generating set of $K(S_{n})$,
confirming Proposition \ref{P-Basis}.

We proceed to describe a class of irreps of $GL_{n}$, that in a formal sense
behave as if they also form $K(S_{n})$.

\section{\textbf{Spherical Principal Series Representations of }$GL_{n}$%
\label{S-SPS}}

In this section we want first to construct/classify the spherical principal
series representations, and second, to recast certain properties of this
collection. Both tasks involve, as in the case of $S_{n}$, the dominance
relation on the set $\mathcal{Y}_{n}$, of Young diagrams with $n$ boxes.

\subsection{\textbf{The Spherical Principal Series}}

Inside $GL_{n}=GL_{n}(\mathbb{F}_{q})$, consider the Borel subgroup $B$ \cite%
{Borel69} of upper triangular matrices%
\begin{equation*}
B=%
\begin{pmatrix}
\ast & \ldots & \ast \\ 
& \ddots & \vdots \\ 
&  & \ast%
\end{pmatrix}%
.
\end{equation*}

Recall, see Section \ref{S-TPP}, that by definition an irreducible
representation $\rho $ of $GL_{n}$ belongs to the \underline{spherical
principal series} (SPS) if it appears inside the induced representation $%
Ind_{B}^{GL_{n}}(\mathbf{1})$, where $\mathbf{1}$ denotes the trivial
representation of $B$.

The construction of the SPS representations, and the verification of some of
their properties can be done intrinsically (e.g., see in Section 10.5. of 
\cite{Gurevich-Howe17}), without the relation to the representation theory
of $S_{n}$. However, for purposes of this note, we prefer to get all the
information from what was obtained already for $S_{n}$ in Section \ref%
{S-RepSn}. This, in particular, will enable us to derive the Pieri rule for $%
GL_{n}$ from that of $S_{n}$.

\subsection{\textbf{The Grothendieck Group of the Spherical Principal Series 
\label{S-GG-SPS}}}

Let us denote by $K_{B}(GL_{n})$ the Abelian group generated, using the
operation of direct sum $\oplus $, from the SPS representations. The notion
of subrepresentation induces a partial order $<$ on $K_{B}(GL_{n})$ and the
intertwining number pairing $\left\langle \text{ },\right\rangle $ (\ref{IN}%
) gives on it an inner product structure.

We proceed to give an effective description of $K_{B}(GL_{n})$.

Recall, see Section \ref{S-TPP}, that the group $K_{B}(GL_{n})$ has a
distinguished collection of members in the form of induced representations
that are associated to Young diagrams. Indeed, to a Young diagram $D\in 
\mathcal{Y}_{n}$ one attaches in a natural a way a flag $F_{D}$ in $\mathbb{F%
}_{q}^{n},$ see Equation (\ref{FD}), and a corresponding parabolic subgroup $%
P_{D}=Stab_{GL_{n}}(P_{D})\subset GL_{n}$. Then, we can consider the trivial
representation $\mathbf{1}$ of $P_{D},$ and induce to obtain%
\begin{equation}
I_{D}=Ind_{P_{D}}^{GL_{n}}(\mathbf{1).}  \label{ID}
\end{equation}

Of course each $I_{D}$ sits inside $Ind_{B}^{GL_{n}}(\mathbf{1})$, but we
can say much more on the relation between the various $I_{D}$'s. Indeed, for
a given Young diagram $D=(d_{1}\geq ...\geq d_{r})\in \mathcal{Y}_{n},$ we
have defined in (\ref{SD}) the subgroup $S_{D}\simeq S_{d_{1}}\times
...\times S_{d_{r}}\subset S_{n}$ and the corresponding Young module $%
Y_{D}=Ind_{S_{D}}^{S_{n}}(\mathbf{1)}$. Then, the Bruhat decomposition \cite%
{Borel69, Bruhat56} gives a bijection between the double cosets%
\begin{equation}
P_{D}\diagdown GL_{n}\diagup P_{E}\text{ \ and \ }S_{D}\diagdown
S_{n}\diagup S_{E},  \label{BD}
\end{equation}%
for every $D,E\in \mathcal{Y}_{n}$.

But, the cardinalities of the double cosets in (\ref{BD}) are exactly the
dimensions of, respectively, the intertwining spaces $%
Hom_{GL_{n}}(I_{D},I_{E})$ and $Hom_{S_{n}}(Y_{D},Y_{E}),$ so we conclude
that,

\begin{proposition}[Bruhat decomposition]
\label{P-Id-IN}For any two Young diagrams $D,E\in \mathcal{Y}_{n}$, we have, 
\begin{equation}
\left\langle I_{D},I_{E}\right\rangle =\left\langle Y_{D},Y_{E}\right\rangle
.  \label{Id-IN}
\end{equation}
\end{proposition}

One way to interpret identity (\ref{Id-IN}) is as follows:

\begin{corollary}
\label{C-i}The correspondence%
\begin{equation}
Y_{D}\longmapsto I_{D},\text{ \ }D\in \mathcal{Y}_{n},  \label{YD-ID}
\end{equation}%
induces an order preserving isometry%
\begin{equation}
\iota :K(S_{n})\widetilde{\rightarrow }K_{B}(GL_{n})\text{.}  \label{iuta}
\end{equation}
\end{corollary}

On how to deduce Corollary \ref{C-i} from Proposition \ref{P-Id-IN}, see the
next section.

\subsection{\textbf{The Grothendieck Groups of }$S_{n}$\textbf{\ and of the
Spherical Principal Series }}

We confirm Corollary \ref{C-i}, and along the way construct the SPS
representations, and deduce various other facts on this collection.

Consider the map $\iota $ (\ref{iuta}), extended by (integral) linearity
from the correspondence (\ref{YD-ID}). Denote by 
\begin{equation}
\rho _{D}=\iota (\sigma _{D}),\text{ \ }D\in \mathcal{Y}_{n}\text{,}
\label{rhoD}
\end{equation}%
the element of $K_{B}(GL_{n})$ corresponding to the irrep (\ref{sigmaD}) of $%
S_{n}$. Note that,

\begin{itemize}
\item $\rho _{D}<I_{D}$; and,

\item $\left\langle \rho _{D},I_{D}\right\rangle =1$,
\end{itemize}

so, in particular, $\rho _{D}$ is irreducible. In fact, the corresponding
properties for $S_{n}$ imply that

\begin{itemize}
\item $\left\langle \rho _{D},Ind_{B}^{GL_{n}}(\mathbf{1})\right\rangle
=\dim (\sigma _{D})$; and,

\item we have, 
\begin{equation}
\{\rho _{D}\}=\widehat{I}_{D}\smallsetminus \bigcup\limits_{D\precneqq E}%
\widehat{I}_{E},  \label{CP-rhoD}
\end{equation}%
i.e., $\rho _{D}$ is the unique irrep that sits in $I_{D}$ (we denote by $%
\widehat{I}_{D}$ the set of irreps inside $I_{D}$) but not in $I_{E}$, for
any Young diagram $E\in \mathcal{Y}_{n}$ that strictly dominates $D$.
\end{itemize}

\begin{remark}
In fact, Property (\ref{CP-rhoD}) characterizes the representation $\rho
_{D} $, and is useful, e.g., you can compute out of it explicitly the
dimension of $\rho _{D}$ and find that (we use bold-face letters to denote
the corresponding algebraic groups \cite{Borel69}) it is equal to $\dim
(\rho _{D})=q^{\dim (\mathbf{GL_{n}/P}_{D})}+o(...),$ as $q\rightarrow
\infty $, a fact that in turn characterizes (again, asymptotically) $\rho
_{D}$ uniquely among all irreps in $I_{D}$.
\end{remark}

How do we know we get all the SPS?\smallskip

A possible answer is that, as we already said, each $I_{D}$ has a unique
irrep that does not occur in the induced module $I_{E}$ corresponding to any
strictly dominating diagram $E\succneqq D$, namely, $\rho _{D}=\iota (\sigma
_{D})$. On the $S_{n}$ side, the irreps $\sigma _{D}$, $D\in \mathcal{Y}_{n}$%
, completely decompose each of the induced representations. By Bruhat, this
transfers to $GL_{n}$, so we get complete decompositions over there also. In
particular, we get a complete decomposition of $Ind_{B}^{GL_{n}}(\mathbf{1}%
)=I_{(1,...,1)}$, the constituents of which are exactly the SPS
representations.

Finally, the above discussion also validates Corollary \ref{C-i}.

Having at our disposal the understanding that the SPS representations and
the representations of $S_{n}$ are in some formal sense the same thing, we
can proceed to discuss the Pieri rule.

\section{\textbf{The Pieri Rule\label{S-PR}}}

Fix $0\leq k\leq n$, and denote by $P_{k,n-k}\subset GL_{n}$ the parabolic
subgroup fixing the first $k$ coordinate subspace of $\mathbb{F}_{q}^{n}$.
There is a natural surjective map $P_{k,n-k}\twoheadrightarrow GL_{k}\times
GL_{n-k}$. Take a Young diagram $D\in \mathcal{Y}_{k}$, and consider the
irreducible SPS representation $\rho _{D}$ of $GL_{k}$ defined by (\ref{rhoD}%
). Denote by $\mathbf{1}_{n-k}$ the trivial representation of $GL_{n-k}$.
Pull back the representation $\rho _{D}\otimes \mathbf{1}_{n-k}$ from $%
GL_{k}\times GL_{n-k}$ to $P_{k,n-k}$ and form the induced representation 
\begin{equation}
I_{\rho _{D}}=Ind_{P_{k,n-k}}^{GL_{n}}(\rho _{D}\otimes \mathbf{1}_{n-k}).
\label{I_rhoD}
\end{equation}

Recall (see Problem \ref{P-PP} in Section \ref{S-TPP}) that, the narrative
of the story we are telling in this note is that, we are seeking to compute
the decomposition of $I_{\rho _{D}}$ (\ref{I_rhoD}) into irreps. Moreover,
it is easy to see that all constituents of the representation $I_{\rho _{D}}$
are SPS, so we are seeking an answer to the decomposition problem in term of
Young diagrams.

To arrive at our goal, after introducing some needed terminology, we
will\smallskip

\textbf{(a) }State the Pieri rule for the long established case of the
complex general linear group $GL_{n}(%
\mathbb{C}
)$.\smallskip

\textbf{(b) }Recall Schur duality.\smallskip

\textbf{(c) }State and prove the Pieri rule for representations of $S_{n}.$

Our proof was suggested by a remark of \textbf{Nolan Wallach}, and uses the
Schur (a.k.a. Schur-Weyl) duality, to deduce the result from the Pieri rule
for $GL_{n}(%
\mathbb{C}
)$.

We note that, recently, in \cite{Ceccherini-Silberstein-Scarabotti-Tolli10},
the authors gave a different proof of the Pieri rule for $S_{n}$. Their
treatment uses the Okounkov-Vershik approach \cite{Okounkov-Vershik05}
through the branching rule from $S_{n}$ to $S_{n-1}$ and Gelfand-Tsetlin
basis.

We would also like to remark that, nowadays the Pieri rule for $S_{n}$ can
be understood as a particular case of the celebrated Littlewood-Richardson
rule \cite{Littlewood-Richardson34, Macdonald79}, but was known \cite%
{Pieri1893}\textit{\ }a long time before this general result.\smallskip
\smallskip

\textbf{(d) }Derive the Pieri rule for the Group $GL_{n}=GL_{n}(\mathbb{F}%
_{q})$, using the equivalence, discussed in Section \ref{S-GG-SPS}, between
the representation theory of the spherical principal series, and that of $%
S_{n}$.

\subsection{\textbf{Skew-Diagrams and Horizontal Strips}}

The various formulations we present of the Pieri rule use the notions of
skew diagram and horizontal strip, that we recall here.

\begin{figure}[h]\centering
\includegraphics
{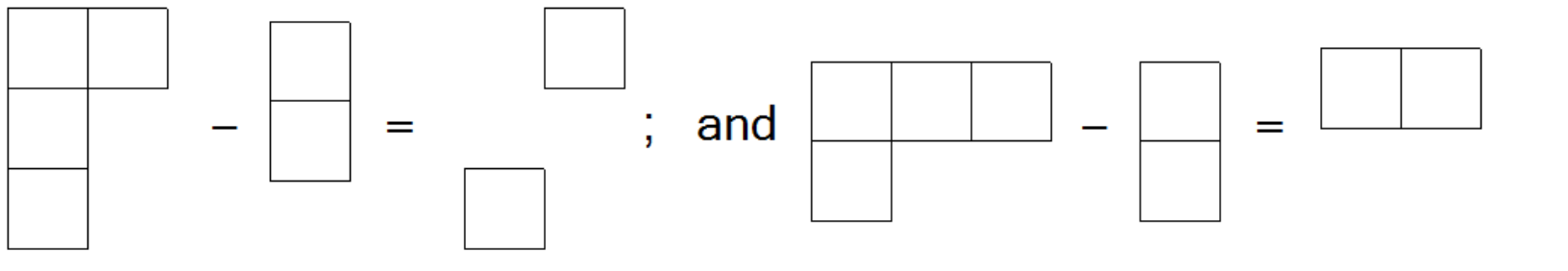}%
\caption{In $\mathcal{Y}_{4}$: $(2,1,1)$, $(3,1),$ contain $(1,1)\in 
\mathcal{Y}_{2},$ with difference a horizontal $2$-strip. }\label{sr}%
\end{figure}

Suppose we have Young diagrams $E\in \mathcal{Y}_{n}$ and $D\in \mathcal{Y}%
_{k}$ such that $E$ contains $D$, denoted $E\supset D$, i.e., each row of $E$
is at least as long as the corresponding row of $D$. Then, by removing from $%
E$ all the boxes belonging to $D$, we obtain a configuration, denoted $E-D,$
called \textit{skew-diagram } \cite{Macdonald79}. If, in addition---see
Figure \ref{sr} for illustration, each column of\ $E$ is at most one box
longer than the corresponding column of $D$, then we call $E-D$ a \textit{%
horizontal strip} (or \textit{horizontal }$m$-\textit{strip }if $E-D$ has $m$
boxes). 

\begin{remark}
In \cite{Ceccherini-Silberstein-Scarabotti-Tolli10} the term that is being
used for "horizontal strip" is "totally disconnected skew diagram".
\end{remark}

\subsection{\textbf{The Pieri Rule for }$GL_{n}(%
\mathbb{C}
)$}

The Pieri rule for $GL_{n}(%
\mathbb{C}
)$ is a (very) special case of the general Littlewood-Richardson rule \cite%
{Howe-Lee12, Littlewood-Richardson34, Macdonald79} for decomposing the
tensor product of any pair of irreducible finite dimensional representations
of $GL_{n}(%
\mathbb{C}
)$. The Pieri rule has been known since the 19th century \cite{Pieri1893},
and is relatively easy to establish \cite{Fulton-Harris91, Howe92, Weyman89}.

There is a standard way to label the irreducible representations of $GL_{n}(%
\mathbb{C}
)$. It is by their highest weights (see, for example \cite{Fulton-Harris91,
Howe92, Weyl46}). A highest weight for $GL_{n}(%
\mathbb{C}
)$ is specified by a decreasing sequence

\begin{equation*}
d_{1}\geq ...\geq d_{n},
\end{equation*}%
of integers.

When all the $d_{j}$ are non-negative, the above sequence can be thought of
as specifying a Young diagram $D$, with $j$-th row having length $d_{j}$.
The number of boxes in $D$ can be arbitrarily large, but the number of rows
is bounded by $n$. Irreps of $GL_{n}(%
\mathbb{C}
)$ corresponding to sequences with all $d_{j}$ non-negative are called 
\textit{polynomial} representations. These are exactly all the irreps of $%
GL_{n}(%
\mathbb{C}
)$ that appear in the tensor powers $(%
\mathbb{C}
^{n})^{\otimes ^{l}}$ of $%
\mathbb{C}
^{n}$ for some $l\geq 0$. (Any irrep of $GL_{n}(%
\mathbb{C}
)$ is isomorphic to a twist by a power of determinant of a polynomial
representation.). We will denote by%
\begin{equation}
\pi _{n}^{D},\text{ \ }D=(d_{1}\geq ...\geq d_{n}\geq 0),  \label{pi-n-D}
\end{equation}%
the polynomial representation of $GL_{n}(%
\mathbb{C}
)$ whose highest weight corresponds to the diagram $D$. The one-rowed
diagrams, given by $(d_{1},0,...,0)$ correspond to the symmetric powers $%
S^{d_{1}}(%
\mathbb{C}
^{n})$.

The Pieri rule for $GL_{n}(%
\mathbb{C}
)$ describes the decomposition of a tensor product $\pi _{n}^{D}\otimes
S^{d}(%
\mathbb{C}
^{n})$ of a general polynomial irrep with a symmetric power $S^{d}(%
\mathbb{C}
^{n})$.

\begin{proposition}[Pieri rule for $GL_{n}(%
\mathbb{C}
)$]
\label{P-PR-GLnC}The representation $\pi _{n}^{D}\otimes S^{d}(%
\mathbb{C}
^{n})$ is multiplicity free. Moreover, we have, 
\begin{equation}
\pi _{n}^{D}\otimes S^{d}(%
\mathbb{C}
^{n})\simeq \sum_{E}\pi _{n}^{E},  \label{PR-GLnC}
\end{equation}%
where $E$ runs through all diagrams such that

\begin{enumerate}
\item $D\subset E;$ and,

\item $E-D$ is an horizontal $d$-strip.
\end{enumerate}
\end{proposition}

\subsection{\textbf{Schur-Weyl Duality}}

The group $GL_{n}(%
\mathbb{C}
)$ is defined in terms of its action on $%
\mathbb{C}
^{n}$. By taking tensor products, this action gives rise naturally to an
action on the $l$-fold tensor product $(%
\mathbb{C}
^{n})^{\otimes ^{l}}$ of $%
\mathbb{C}
^{n}$ with itself (a.k.a., the $l$-th tensor power of $%
\mathbb{C}
^{n}$). Clearly, the permutation group $S_{l}$ also acts on $(%
\mathbb{C}
^{n})^{\otimes ^{l}}$ by permuting the factors of the product. This action
of $S_{l}$ clearly commutes with the action of $GL_{n}(%
\mathbb{C}
)$. Schur-Weyl duality \cite{Howe92, Schur27, Weyl46} says that

\begin{proposition}[Schur-Weyl duality - non-explicit form]
The actions of $S_{l}$ and $GL_{n}(%
\mathbb{C}
)$ on $(%
\mathbb{C}
^{n})^{\otimes ^{l}}$generate mutual commutants of each other.
\end{proposition}

From this, Burnside's double commutant theorem \cite{Burnside1905, Weyl46}
lets us conclude that, as an $S_{l}\times $ $GL_{n}(%
\mathbb{C}
)$-module, we have a decomposition

\begin{equation}
(%
\mathbb{C}
^{n})^{\otimes ^{l}}\simeq \sum_{D}\sigma _{l}^{D}\otimes \tau _{n}^{D},
\label{SWD-1}
\end{equation}%
where $D\in \mathcal{Y}_{l}$ runs through diagrams with $l$ boxes, $\sigma
_{l}^{D}$ are the associated irreps (\ref{sigmaD}) of $S_{l},$ and the $\tau
_{n}^{D}$ are appropriate irreps of $GL_{n}(%
\mathbb{C}
)$. Some computation then shows that, remarkably, $\tau _{n}^{D}$ is equal
to the representation $\pi _{n}^{D}$ (provided of course that $D$ does not
have more than $n$ rows; otherwise $\tau _{n}^{D}=0$) given by Equation (\ref%
{pi-n-D}). Thus, we can rewrite (\ref{SWD-1}), and obtain

\begin{proposition}[Schur-Weyl duality - explicit form]
\label{P-SWD-EF}As an $S_{l}\times $ $GL_{n}(%
\mathbb{C}
)$-module, we have the decomposition%
\begin{equation}
(%
\mathbb{C}
^{n})^{\otimes ^{l}}\simeq \sum_{D}\sigma _{l}^{D}\otimes \pi _{n}^{D},
\label{SWD-2}
\end{equation}%
where $D$ runs over all diagrams in $\mathcal{Y}_{l}$ with at most $n$ rows.
\end{proposition}

\subsection{\textbf{The Pieri Rule for }$S_{n}$}

With the usual notation, consider $k<n$, and $S_{k}\subset S_{n}$, in the
standard way, as the group that fixes the last $n-k$ letters on which $S_{n}$
acts. Then the symmetric group on these letters is $S_{n-k}$, and we have
the product $S_{k}\times S_{n-k}\subset S_{n}.$

Take a partition/Young diagram $D$ of size $k$, and let $\sigma _{D}$ be the
associated irreducible representation (\ref{sigmaD}) of $S_{k}$. Let $%
\mathbf{1}_{n-k}$ be the trivial representation of $S_{n}$. Form the induced
representation%
\begin{equation}
I_{\sigma _{D}}=Ind_{S_{k}\times S_{n-k}}^{S_{n}}(\sigma _{D}\otimes \mathbf{%
1}_{n-k}),  \label{I_sigmaD}
\end{equation}%
of $S_{n}$.

The Pieri Rule for $S_{n}$ describes the decomposition of this induced
representation into irreducible subrepresentations.

\begin{theorem}[Pieri rule for $S_{n}$]
\label{T-PR-Sn}The representation $I_{\sigma _{D}}$ (\ref{I_sigmaD}) is
multiplicity-free. It consists of one copy of each representation $\sigma
_{E}$ of $S_{n}$, for diagrams $E\in \mathcal{Y}_{n},$ such that

\begin{enumerate}
\item $D\subset E;$ and,

\item $E-D$ is an horizontal $(n-k)$-strip.
\end{enumerate}
\end{theorem}

In Appendix \ref{P-T-PR-Sn} we give our proof of Theorem \ref{T-PR-Sn},
demonstrating how it follows from the Pieri rule for $GL_{n}(%
\mathbb{C}
)$, invoking the Schur-Weyl duality.

\begin{remark}[Description of the Young Module]
\label{R-D-YD}Theorem \ref{T-PR-Sn} can be used to give a recursive
description of the Young module $Y_{D}$ (\ref{YD}).

Given a Young diagram $D\in \mathcal{Y}_{n}$ with $n$ boxes, let $D_{s}$ be
the diagram consisting of the first $s$ rows of $D$, and let $k_{s}$ be the
number of boxes in $D_{s}$. Suppose that in $D$ there are $r$ rows in all,
so that $k_{r}=n$. Suppose we know how to decompose $Y_{D_{s}}$. Then, if we
apply the Pieri rule to each component of $Y_{D_{s}}$ and $k_{s+1}$ is the
number of boxes in $D_{s+1}$, we learn how to decompose $Y_{D_{s+1}}$.
Starting with $s=1$, we can successively decompose the $Y_{D_{s}}$ for all $%
s $ up to $r$, at which point we will have found the decomposition of $Y_{D}$%
.

For example, the above method provides us with the following combinatorial
description of the multiplicity of the irrep $\sigma _{E},$ $E\in \mathcal{Y}%
_{n}$, in $Y_{D}$: it is the number of ways to fill $E$ with a nested family
of sub-diagrams $E_{s}$, such that

\begin{itemize}
\item $E_{s}\subset E_{s+1}$; and,

\item $E_{s+1}-E_{s}$ is a horizontal strip with $k_{s+1}-k_{s}$ boxes.
\end{itemize}

From this, we can see, again, that the multiplicity of $\sigma _{D}$ in $%
Y_{D}$ is $1$.
\end{remark}

\subsection{\textbf{The Pieri Rule for }$GL_{n}(\mathbb{F}_{q})$}

Now we can finish our story, and deliver the answer to the introduction's
motivating problem of decomposing $I_{\rho _{D}}$ (\ref{I_rhoD}).

Note that the isomorphism $\iota $ (\ref{iuta}), between the representation
groups $K(S_{n})$ and $K_{B}(GL_{n}),$ sends $I_{\sigma _{D}}$ to $I_{\rho
_{D}}$. So the Pieri rule for $S_{n},$ implies the Pieri rule for $GL_{n}(%
\mathbb{F}_{q})$, i.e., the same description as in Theorem \ref{T-PR-Sn},
just replace there, $S_{n}$ by $GL_{n}(\mathbb{F}_{q}),$ and $\sigma _{D},$ $%
\sigma _{E},$ by $\rho _{D},\rho _{E},$ respectively.

\begin{remark}[Decomposing permutation repesentation on flag variety]
Replacing the Young module $Y_{D}$, in Remark \ref{R-D-YD}, by the induced
representation $I_{D}=Ind_{P_{D}}^{GL_{n}}(\mathbf{1),}$ which is the the
space of functions on the flag variety $GL_{n}/P_{D}.$ We get a recursive
formula for the decomposition into irreps of the permutation representation
of $GL_{n}$ on functions on a fairly general flag variety.
\end{remark}

\appendix

\section{\textbf{Proofs}}

\subsection{\textbf{Proofs for Section \protect\ref{S-RepSn}}}

\subsubsection{\textbf{Proof of Theorem \protect\ref{T-Main}\label{P-T-Main}}%
}

For a set $X$ let us denote by $L(X)$ the space of complex valued functions
on $X$. We also use this notation to denote the standard permutation
representation of a group $G$, in case it acts on $X$.

Now we can proceed to give the proof.

\begin{proof}
\textbf{Part (1).}\textit{\ }Let us analyze the space of intertwiners $%
Hom(Y_{E}(sgn),Y_{D})$. This has a "geometric" description from which the
information we are after can be read.

\begin{itemize}
\item First, recall that we can realize $Y_{D}$ as the permutation
representation $L(\mathcal{T}_{D})$ associated with the action of $S_{n}$ on
the set $\mathcal{T}_{D}$ of all tabloids that one can make out of $D$ (see
Section \ref{S-YM}). In the same way, $Y_{E}(sgn)$ can be realized on the
space $L(\mathcal{T}_{E})$ with the permutation action of $S_{n}$ on it
twisted by $sgn$.

\item Second, using the bases of delta functions of $L(\mathcal{T}_{E})$ and 
$L(\mathcal{T}_{D}),$ we can associate to every intertwiner $%
Hom(Y_{E}(sgn),Y_{D})$ a kernel function (i.e., a matrix) $K$ on $\mathcal{T}%
_{D}\times \mathcal{T}_{E}$ that satisfies 
\begin{equation}
K(s(T_{D}),s(T_{E}))=sgn(s)K(T_{D},T_{E}),  \label{K}
\end{equation}%
for every $s\in S_{n}$, $T_{D}\in \mathcal{T}_{D}$ and $T_{E}\in \mathcal{T}%
_{E}$.

Let us denote by $L(\mathcal{T}_{D}\times \mathcal{T}_{E})^{\mathbf{1}%
\otimes sgn}$ the collection of all $K$ satisfying Identity (\ref{K}).
\end{itemize}

In summary, we obtained, 
\begin{equation}
Hom(Y_{E}(sgn),Y_{D})=L(\mathcal{T}_{D}\times \mathcal{T}_{E})^{\mathbf{1}%
\otimes sgn}\text{.}  \label{K-space}
\end{equation}%
This is a geometric description of the space of intertwiners.

Now suppose we have $K$ from (\ref{K-space}), and suppose there are $%
T_{D}\in \mathcal{T}_{D}$ and $T_{E}\in \mathcal{T}_{E}$ with rows, one of $%
T_{D}$ and one of $T_{E}$, that share two numbers $i,j\in \{1,...,n\}$.
Then, the permutation that transposes $i$ and $j$ must preserve $%
K(T_{D},T_{E})$, and also must change its sign. Therefore, $%
K(T_{D},T_{E})=-K(T_{D},T_{E})$, so $K(T_{D},T_{E})=0$. In other words $%
K(T_{D},T_{E})\neq 0$ only if

\begin{itemize}
\item each number from the first row of $E,$ should sit in a different row
of $D,$ so, \smallskip

$E_{1}\leq D_{1}^{t},\smallskip $

i.e., the length $E_{1}$ of the first row of $E$ is not more than that of
the first row of $D^{t}$.\smallskip\ 

and\smallskip

\item each number from the second row of $E$ should sit in a different row
of $D,$ so we also have,\smallskip

$E_{1}+E_{2}\leq D_{1}^{t}+D_{2}^{t}.\smallskip $

\item etc...
\end{itemize}

Namely, for the space (\ref{K-space}) to be non-trivial, it is necessary to
have\smallskip\ 

$(\ast )$ $\ E\preceq D^{t}$.\medskip

Next, assuming $E=D^{t},$ we want to show that the intertwining space (\ref%
{K-space}) is one dimensional.

Let us first give one orbit in $\mathcal{T}_{D}\times \mathcal{T}_{D^{t}}$
that supports a non-trivial intertwiner:

\begin{itemize}
\item Take the Young diagram $D$ and fill each box of it with numbers from $%
\{1,..,n\}$. The object we obtained in this way is called Young \textit{%
tableau }\cite{Fulton97}. From it we can make in a natural way a Young
tabloid $T_{D}$ by grouping together the numbers in each line of the tableau.

\item We could also first "transpose" the filled $D$ to obtain a tableau
associated with $D^{t},$ and then, in the same way as above, form the
corresponding tabloid $T_{D^{t}}$.
\end{itemize}

It is clear that, any two rows, one of $T_{D}$ and one of $T_{D^{t}},$ share
no more than one number in common. Hence, the group $S_{n}$ acts freely on
the orbit 
\begin{equation}
\mathcal{O}_{T_{D},T_{D^{t}}}\subset \mathcal{T}_{D}\times \mathcal{T}%
_{D^{t}}  \label{O-TDt-TD}
\end{equation}%
of $(T_{D},T_{D^{t}})$, and, in particular, there exists an intertwiner $K$
from (\ref{K-space}) which is supported on it.

Now, let us show that (\ref{O-TDt-TD}) is the only orbit that supports such $%
K$. Indeed, take $s\in S_{n}$, such that $s(T_{D^{t}})\neq T_{D^{t}}$. But
then, there are rows, one of $s(T_{D^{t}}),$ and one of $T_{D},$ that share
two numbers in common, then, as was explained earlier, the orbit of $%
(T_{D},s(T_{D^{t}}))$ does not support an intertwiner.\medskip

Finally, let us show that if $E\preceq D^{t}$, then the space (\ref{K-space}%
) is non-zero, i.e., the condition $(\ast )$ is also sufficient. It is
enough to examine the case when $E=(D^{t})^{\circ }$ obtained from $D^{t}$,
by moving one box down to form a new lower row. Now, look at the tabloids $%
T_{D}$ and $T_{D^{t}},$ that we used in the paragraph just above, and the
natural tabloid $T_{(D^{t})^{\circ }}$ one obtains from the filling of $%
D^{t} $ by numbers as we did above in order to create $T_{D^{t}}$. Then, as
we argued above, the orbit of $(T_{D},T_{(D^{t})^{\circ }})$ supports a
non-trivial intertwiner.\medskip

\textbf{Part (2). }If $D$ is not dominated by $E$, then from Part (1) we see
that $\left\langle Y_{E^{t}}(sgn),Y_{D}\right\rangle =0$, and in particular,
again by Part (1), $Y_{E}$ cannot be a subrepresentation of $Y_{D}$%
.\smallskip\ 

On the other hand, let us assume that $D$ is strictly dominated by $E$, and
show that $Y_{E}\lvertneqq Y_{D}$.

First, we realize the space of intertwiners between $Y_{E}$ and $Y_{D}$
geometrically,%
\begin{equation*}
Hom(Y_{E},Y_{D})=L(\mathcal{T}_{D}\times \mathcal{T}_{E})^{S_{n}},
\end{equation*}%
where on the right hand side of the equality we have, the space of $S_{n}$%
-invariant kernels $K$ on $\mathcal{T}_{D}\times \mathcal{T}_{E}$, or
equivalently the space of functions on the set of orbits $S_{n}\diagdown (%
\mathcal{T}_{D}\times \mathcal{T}_{E}).$

Second, we can parametrize the above set of orbits as follows. Take $%
T_{D}\in \mathcal{T}_{D},$ and $T_{E}\in \mathcal{T}_{E},$ and denote by $%
R_{i}(T_{D})$ and $R_{j}(T_{E}),$ the $i$-th row of $T_{D}$, and $j$-th row
of $T_{E},$ respectively. Then, we can define the \textit{intersection matrix%
} 
\begin{equation}
R_{T_{D},T_{E}}=(r_{ij}),\text{ \ \ }r_{ij}=\#\left( \text{ }%
R_{i}(T_{D})\cap R_{j}(T_{E})\right) ,  \label{R-TD-TE}
\end{equation}%
i.e., $r_{ij}$ is the number of elements common to both rows. It is clear
that $R_{T_{D},T_{E}}$ is an invariant of the orbit. Moreover, it gives a
complete invariant. Indeed, it is not difficult to see that if $%
R_{T_{D},T_{E}}=R_{T_{D}^{\prime },T_{E}^{\prime }}$, then there exists $%
s\in S_{n}$ such that $s(T_{D})=T_{D}^{\prime },$ and $s(T_{E})=T_{E}^{%
\prime }$.

A direct computation, using the parametrization (\ref{R-TD-TE}), reveals
that,

\begin{claim}
Consider the Young diagrams $D_{n-k,k}=(n-k,k)$ and $D_{n-k^{\prime
},k^{\prime }}=(n-k^{\prime },k^{\prime })$, where $0\leq k,k^{\prime }\leq 
\frac{n}{2}$. Then, 
\begin{equation*}
\left\langle Y_{D_{n-k,k}},Y_{D_{n-k^{\prime },k^{\prime }}}\right\rangle
=\min \{k+1,k^{\prime }+1\}\text{.}
\end{equation*}
\end{claim}

So $Y_{D_{n,0}}$ contains $1$ representation - the trivial representation.
Then $Y_{D_{n-1,1}}$ contains two representations, one of which is the
trivial representation. Since $Y_{D_{n-2,2}}$ has intertwining number $1$
with $Y_{D_{n,0}}$ and $2$ with $Y_{D_{n-1,1}}$, it must contain the two
representations of $Y_{D_{n-1,1}}$ with multiplicity $1$ each. Since its
self intertwining number is $3$, it contains $3$ representations, each with
multiplicity $1$. Then we can continue like that: $Y_{D_{n-3,3}}$ contains
each of the representations of $Y_{D_{n-2,2}}$ with multiplicity $1$, and
then one new representation, and so on. So in particular:\smallskip\ 

$(\ast \ast )$ $\ Y_{D_{n-k-1,k+1}}$ contains $Y_{D_{n-k,k}}$ when $k+1\leq 
\frac{n}{2}.\smallskip $

Now take any diagram $D$, containing two rows $R$ and $R^{\prime }$, with $%
R^{\prime }$ (which might be of length equal to $0$) at least two boxes
shorter than $R$. Then we can form $Y_{D}$ by first forming the
representation $Y_{D_{R,R^{\prime }}}$ of $S_{R+R^{\prime }}$, and then
extending to be trivial on the stabilizers of the other rows, and then
inducing up to $S_{n}$. So if we replace $R$ and $R^{\prime }$ with $R-1$
and $R^{\prime }+1$, we will get a larger representation, using Fact $(\ast
\ast ).$ This completes the verification of Part (2), and of Theorem \ref%
{T-Main}.
\end{proof}

\subsubsection{\textbf{Proof of Corollary \protect\ref{C-Class}\label%
{P-C-Class}}}

\begin{proof}
Note that the dominance order on $\mathcal{Y}_{n}$ is a partial order, and
in particular, is anti-symmetric, i.e., for every $E,D\in \mathcal{Y}_{n}$,
if $E\preceq D$ and $E\succeq D$, then $E=D$. But, if $\sigma _{E}\simeq
\sigma _{D}$, then, by the "iff" of Part (1) of Theorem \ref{T-Main}, $%
E\preceq D$ and $E\succeq D,$ so the Corollary follows.
\end{proof}

\subsection{\textbf{Proofs for Section \protect\ref{S-PR}}}

\subsubsection{\textbf{Proof of Theorem \protect\ref{T-PR-Sn}\label%
{P-T-PR-Sn}}}

\begin{proof}
Schur duality for $S_{k}\times GL_{n}(%
\mathbb{C}
)$ on the $k$-fold tensor product $(%
\mathbb{C}
^{n})^{\otimes ^{k}}$ says (Proposition \ref{P-SWD-EF}, Equation (\ref{SWD-2}%
)) that we have 
\begin{equation*}
(%
\mathbb{C}
^{n})^{\otimes ^{k}}\simeq \sum_{D\in \mathcal{Y}_{k}}\sigma _{k}^{D}\otimes
\pi _{n}^{D},
\end{equation*}%
where $\pi _{n}^{D}$ is the irrep of $GL_{n}(%
\mathbb{C}
)$ with highest weight corresponding to the diagram $D$.

We can also apply Schur duality to the action of $S_{n-k}\times GL_{n}(%
\mathbb{C}
)$ on $(%
\mathbb{C}
^{n})^{\otimes ^{n-k}}$. Then the space of fixed vectors for $S_{n-k}$ is
the $S_{n-k}\times GL_{n}(%
\mathbb{C}
)$-module $\mathbf{1}_{n-k}\otimes \pi _{n}^{(n-k)}$, corresponding to the
diagram with one row of length $n-k$. This is just the $(n-k)$-th symmetric
power of the standard action on $%
\mathbb{C}
^{n}$.

Now consider,%
\begin{equation*}
(%
\mathbb{C}
^{n})^{\otimes ^{n}}\simeq (%
\mathbb{C}
^{n})^{\otimes ^{k}}\otimes (%
\mathbb{C}
^{n})^{\otimes ^{n-k}},
\end{equation*}%
as an $S_{n}$-module, again with $S_{n}$ acting by permutation of the
factors. If we take the isotypic component of $\sigma _{D}$ inside $(%
\mathbb{C}
^{n})^{\otimes ^{k}}$, and the space of fixed vectors for $S_{n-k}$ inside $(%
\mathbb{C}
^{n})^{\otimes ^{n-k}}$, their tensor product will be the isotypic component
inside $(%
\mathbb{C}
^{n})^{\otimes ^{n}}$ of the representation $\sigma _{D}\otimes \mathbf{1}%
_{n-k}$ of $S_{k}\times S_{n-k}$.

On the other hand, the action of $GL_{n}(%
\mathbb{C}
)$ on the indicated tensor product is described by (a multiple of) the
tensor product $\pi _{n}^{D}\otimes \pi _{n}^{(n-k)}$. By the Pieri rule for 
$GL_{n}(%
\mathbb{C}
)$ (see Proposition \ref{P-PR-GLnC}) this decomposes into a multiplicity
free sum of irreps for $GL_{n}(%
\mathbb{C}
)$ whose highest weights are given by diagrams $E$ having the form indicated
in the statement of the proposition: $E$ has $n$ boxes, contains $D$, and $%
E-D$ consists of a horizontal $(n-k)$-strip. Thus, the $S_{k}\times S_{n-k}$
isotypic component for $\sigma _{D}\otimes \mathbf{1}_{n-k}$ of $(%
\mathbb{C}
^{n})^{\otimes ^{n}}$ has the structure%
\begin{equation}
\sum_{E}\sigma _{D}\otimes \mathbf{1}_{n-k}\otimes \pi _{n}^{E},  \label{Sum}
\end{equation}%
as $S_{k}\times S_{n-k}\times GL_{n}(%
\mathbb{C}
)$-module, where $E$ runs over the diagrams specified in the statement of
the proposition.

Now consider the representation of $S_{n}$ generated by this space. By Schur
duality for $S_{n}$, it will be%
\begin{equation*}
\sum_{E}\sigma _{E}\otimes \pi _{n}^{E},
\end{equation*}%
as $S_{n}\times GL_{n}(%
\mathbb{C}
)$-module. Comparing this with Formula (\ref{Sum}), we conclude that each
representation $\sigma _{E}$ of $S_{n}$ contains one copy of the
representation $\sigma _{D}\otimes \mathbf{1}_{n-k}$ when restricted to $%
S_{k}\times S_{n-k}$, and that these are the only representations of $S_{n}$
that do contain $\sigma _{D}\otimes \mathbf{1}_{n-k}$. By Frobenius
reciprocity, this is equivalent to the statement of Theorem \ref{T-PR-Sn}.
The proof is complete.
\end{proof}


\begin{thebibliography}{Ceccherini-Silberstein-Scarabotti-Tolli10}
\bibitem[Borel69]{Borel69} Borel A., Linear algebraic groups. \textit{GTM
126, Springer-Verlag (1969).}

\bibitem[Bruhat56]{Bruhat56} Bruhat F., Sur les representations induites des
groupes de Lie. \textit{Bull. Soc. Math. France 84 (1956) 97-205.}

\bibitem[Burnside1905]{Burnside1905} Burnside W., On the condition of
reducibility of any group of linear substitutions. \textit{Proc. London
Math.Soc. 3 (1905) 430-434.}

\bibitem[Ceccherini-Silberstein-Scarabotti-Tolli10]%
{Ceccherini-Silberstein-Scarabotti-Tolli10} Ceccherini-Silberstein T.,
Scarabotti F., and Tolli F., Representation Theory of the Symmetric Groups:
The Okounkov-Vershik Approach, Character Formulas, and Partition Algebras. 
\textit{Cambridge Studies in Advanced Mathematics 121 (2010).}

\bibitem[Frobenius68]{Frobenius68} Frobenius F.G., Gesammelte Abhandlungen. 
\textit{Springer-Verlag, Berlin, 1968.}

\bibitem[Fulton97]{Fulton97} Fulton W., Young Tableaux. \textit{Lon. Mat.
Soc. Stu. Tex. 35, Cambridge University Press (1997).}

\bibitem[Fulton-Harris91]{Fulton-Harris91} Fulton W. and Harris J.,
Representation theory: A first course. \textit{GTM 129, Springer (1991).}

\bibitem[Gurevich-Howe17]{Gurevich-Howe17} Gurevich S and Howe R., Rank and
Duality in Representation Theory. \textit{Takagi lectures Vol. 19, Japanese
Journal of Mathematics (2017). }

\bibitem[Gurevich-Howe19]{Gurevich-Howe19} Gurevich S and Howe R., Harmonic
Analysis on $GL_{n}$ over Finite Fields. \textit{Submitted (2019).}

\bibitem[Howe92]{Howe92} Howe R., Perspectives on invariant theory: Schur
duality, multiplicity-free actions and beyond. \textit{The Schur lectures
(1992) 1--182.}

\bibitem[Howe-Moy86]{Howe-Moy86} Howe R. and Moy A., Harish-Chandra
homomorphisms for p-adic groups. \textit{CBMS Regional Conference Series in
Mathematics 59 (1986).}

\bibitem[Howe-Lee12]{Howe-Lee12} Howe, R. and Lee S.T., Why should the
Littlewood--Richardson rule be true. \textit{Bull. Am. Math. Soc. 43,
187--236 (2012).}

\bibitem[Littlewood-Richardson34]{Littlewood-Richardson34} Littlewood D. and
Richardson A., Group characters and algebra. \textit{Philos. Trans. Roy.
Soc. London A 233 (1934) 99-142.}

\bibitem[Macdonald79]{Macdonald79} Macdonald I.M., Symmetric functions and
Hall polynomials. \textit{Oxford Mathematical Monographs (1979).}

\bibitem[Mackey51]{Mackey51} Mackey G.W., On induced representations of
groups. \textit{Amer. J. Math. 73 (1951), 576--592.}

\bibitem[Okounkov-Vershik05]{Okounkov-Vershik05} Okounkov, A.Y and Vershik
A.M., A New Approach to the Representation Theory of the Symmetric Groups.
II. \textit{J Math Sci 131, 5471--5494 (2005).}

\bibitem[Pieri1893]{Pieri1893} Pieri M., Sul problema degli spazi secanti. 
\textit{Rend. Ist. Lombardo 26 (1893) 534--546.}

\bibitem[Sagan91]{Sagan91} Sagan B., The Symmetric Group. \textit{%
Springer-Verlag New York, Inc. (2001).}

\bibitem[Schur27]{Schur27} Schur I., Uber die rationalen Darstellungen der
allgemeinen linearen Gruppe. \textit{Springer (1927).}

\bibitem[Serre77]{Serre77} Serre J.P., Linear Representations of Finite
Groups. \textit{Springer (1977).}

\bibitem[Weyl46]{Weyl46} Weyl H., The classical groups, their invariants and
representations. \textit{Princeton Univ. Press (1946).}

\bibitem[Weyman89]{Weyman89} Weyman J., Pieri's formulas for classical
groups. \textit{Contempory Math 188, American Mathematical Society,
Providence RI (1989) 177-184.}
\end{thebibliography}
\end{document}